		\documentclass[reqno, letterpaper,11pt]{amsart}  
 
		\usepackage{amsmath, amssymb}
		\usepackage[all, knot, frame]{xy}
			\xyoption{arc} 
		\usepackage{tikz}
			\usetikzlibrary {snakes}
		\usepackage{natbib}

		\bibpunct{[}{]}{,}{a}{}{;} 
						
			\def\a{\alpha }
			\def\b{\beta }

			\def\g{\gamma }

			\def\la{\lambda }
			\def\s{\sigma }

			\def\N{{\mathbb N}}
			\def\Z{{\mathbb Z}}
			\def\Q{{\mathbb Q}}

			\def\<{\left\langle}
			\def\>{\right\rangle}
	
			\newtheorem{thm}{Theorem}[section]
			
			\newtheorem{proposition}[thm]{Proposition}

			\theoremstyle{definition}
			\newtheorem{definition}[thm]{Definition}
			
			\newtheorem{example}[thm]{Example}

\begin{document}

\title{Obstructing Sliceness in a Family of Montesinos Knots}
\author{Luke Williams}
\address{Department of Mathematics, Michigan State University, East Lansing, MI 48824}
\email{will2086@msu.edu}

	\begin{abstract}
		Using Gauge theoretical techniques employed by Lisca for 2-bridge knots and by 
		Greene-Jabuka for 3-stranded pretzel knots, we show that no member of the family of 
		Montesinos knots $M(0;[m_1+1,n_1+2],[m_2+1,n_2+2],q)$, with certain restrictions on 
		$m_i$, $n_i$, and $q$, can be (smoothly) slice.  Our techniques use Donaldson's 
		diagonalization theorem and the fact that the 2-fold covers of Montisinos knots bound plumbing 
		4-manifolds, many of which are negative definite.  Some of our examples include knots 
		with signature 0 and square determinant.
	\end{abstract}
	
\maketitle


\section{Introduction}
	In his recent breakthrough articles \cite{Lisca-2007, Lisca-2007-2}, Lisca applies gauge 
	theory, based on work of Donaldson \cite{Donaldson-1987}, to obstruct smooth sliceness for
	2-bridge knots.  Lisca's approach uses the observation that the 2-fold branched cover of
	every 2-bridge knot bounds a negative definite plumbing 4-manifold.  On the other hand, if
	a 2-bridge knot is slice, its 2-fold cover also bounds a rational homology 4-ball.  Gluing 
	these two 4-manifolds along their common boundary yields a smooth, closed, negative definite 
	4-manifold $X$.  As such, according to Donaldson, its intersection form has to be standard.
	As Lisca goes on to show, this obstruction suffices to pin down all slice 2-bridge knots and
	moreover, suffices to determine the smooth concordance orders of all 2-bridge knots.   
	
	In \cite{Greene-Jabuka}, Greene and Jabuka use this approach, supplemented by Heegaard Floer 
	homology techniques, to determine all slice knots among the 3-stranded pretzel knots $P(p,q,r)$
	with $p,q,r$ odd.  Both Lisca's article, in the case of 2-bridge knots, and Greene-Jabuka's article,
	in the case of 3-stranded pretzel knots, resolve in the affirmative the slice-ribbon conjecture 
	\cite{kirby}.
	
	Building on Lisca's work, we employ the same gauge theoretic techniques to address the question
	of smooth sliceness for a five parameter family of Montesinos knots with three rational tangles.  
	Specifically, the family we consider is (see Section \ref{mont def} for a precise definition)
	\[
		M(0;[m_1+1,n_1+2],[m_2+1,n_2+2],q)
	\]
	with $m_i, n_i, -q\in \N$.  Our main theorem is
	\begin{thm} \label{main}
		Let $m_1, n_1, m_2, n_2, -q\in \N$ subject to the condition 
		\[
			\frac{n_1+2}{\left(m_1+1\right)n_1+2m_1+1}+\frac{n_2+2}{\left(m_2+1\right)n_2+2m_2+1}+\frac{1}{q}>0.
		\]
		Moreover, assume that $m_1, n_1, m_2, n_2, q$ were chosen so that $K=M(0;[m_1+1,n_1+2],[m_2+1,n_2+2],q)$
		is a knot (see Proposition \ref{m_i n_i q parity proposition}).  Then $K$ is not smoothly slice.
	\end{thm}
	The conditions on $m_i, n_i, q$ are realized on a multitude of examples.  We only list a few but the 
	interested reader will have little trouble finding many additional examples.
	\begin{example}
	Choose $q=-3$, $m_1=1$, and $m_2=2$.  Then the conditions from Theorem \ref{main} translate into 
	\[
		n_1, n_2\geq 1 \hspace{1cm} \mbox{and} \hspace{1cm} \mbox{$n_1$ is even or $n_2$ is odd}.
	\]
	Consequently there are infinitely many examples among $M(0;[2,n_1+2],[3,n_2+2],-3)$ which are obstructed from being
	smoothly slice by Theorem \ref{main}.
	\end{example}	
	Before stating our next example, we remind the reader that the knot signature $\s(K)$ and knot determinant $\det(K)$
	of a knot $K$ can used as obstructions to sliceness.  If $K$ is slice then $\s(K)=0$ and $|\det(K)|$ is a square.  
	With this in mind, consider the next example.
	\begin{thm}\label{infinite family}
		The family $M(0;[m_1+1,n_1+2],[m_2+1,n_2+2],q)$ (with the restrictions on $m_i$, $n_i$ and $q$ as outlined 
		in Theorem \ref{main}) contains an infinite subfamily of knots whose signatures are 0 and whose determinants 
		are square.
	\end{thm}
	
	\subsection{Organization}
	In Section \ref{Obstruction Section}, we formally outline the approach employed by Lisca to 
	obstruct smooth sliceness.  In Section \ref{mont def}, we give definitions of a Montesinos knot
	as well as outline plumbing 4-manifolds whose boundaries are 2-fold branched covers of such knots.
	Finally, in Sections \ref{main prof section} and \ref{infinite subfamily section}, we prove Theorems \ref{main}
	and \ref{infinite family} respectively.


\section{Obstruction to Sliceness}\label{Obstruction Section}
		This section outlines the technique used by Lisca \cite{Lisca-2007}.  To start with,
		let $K$ be a knot in $S^3$ and let $Y_K$ denote the 3-maniold obtained by taking the
		the double branched cover of $S^3$ with branching set $K$.  We note that $Y_K$ is 
		always a rational homology sphere and moreover, if $K$ is slice\footnote{All instances of 
		sliceness in the article will always refer to smooth sliceness.}, then $Y_K$ bounds
		a rational homology 4-ball $W_K$.  The manifold $W_K$ is obtained by taking a double
		branched cover of the 4-ball branched over the slice disk of $K$.  
		
		If $K$ is a knot whose double branched cover $Y_K$ bounds a negative definite smooth 
		4-manifold $X_K$, one can form a closed, smooth, negative definite 4-manifold by gluing
		$W_K$ to $X_K$ along $Y_K$.  Recall Donaldson's celebrated theorem  
	\begin{thm}[{Donaldson's Theorem A} \cite{Donaldson-1987}]\label{Donaldson's Theorem}
		Let $W$ be a closed, oriented, smooth 4-manifold with negative definite intersection
		form $Q_W$, then $Q_W$  is isomorphic to the standard negative definite form $-I$ of 
		the same rank.
	\end{thm}
		Combining Donaldson's theorem with the remarks preceding it, we arrive at 
	\begin{thm}\label{main obstruction}
		If $K$ is a slice knot whose double branched cover $Y_K$ bounds the negative definite
		smooth 4-manifold $X_K$, then the intersection form of $X_K$ embeds into the standard 
		negative definite form of the same dimension.  
	\end{thm}
	\begin{example}\label{figure eight}
		Consider $K$ the figure-eight knot.  $Y_K$ bounds a plumbing 4-manifold given by a single 
		component of the plumbing given in Figure \ref{figure eight plumbing}.  If $K$ were slice,
		then there would be an embedding (see Section \ref{main proof} for notation)
		$\varphi:\Z^2\to \Z^2$ such that $\varphi(f_1)^2=-3$.  No
		such embedding exists since there are only two basis elements in the codomain.  
		So, Theorem \ref{main obstruction} implies that $K$ could not be slice.  
		
		We know that $K\# K$ is slice.  
		$Y_{K\# K}$ bounds the plumbing 4-manifold in Figure \ref{figure eight plumbing}.  With
		the indicated basis, we can find an embedding $\varphi:\Z^4\to \Z^4$:
		\[
			\varphi(f_1)=e_2+e_3+e_4, 	\hspace{.2in}
			\varphi(f_2)=e_1-e_2, 		\hspace{.2in}
			\varphi(f_3)=e_1+e_2-e_3, 	\hspace{.2in}
			\varphi(f_4)=e_3-e_4.
		\]
		Then, as expected, we cannot use Theorem \ref{main obstruction} to obstruct $K\# K$ from
		being slice.  
		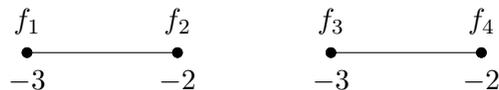
\begin{figure}[htp]
		{
			\[
				\begin{tikzpicture} 
					\draw[xshift=-1in] (-2, .0) -- ( 0, .0); 
					\fill[xshift=-1in] (-2, .0) circle (2pt) node[below =3pt] {$-3$};
					\fill[xshift=-1in] ( 0, .0) circle (2pt) node[below =3pt] {$-2$};
					\fill[xshift=-1in] (-2, .0) circle (2pt) node[above =3pt] {$f_1$};
					\fill[xshift=-1in] ( 0, .0) circle (2pt) node[above =3pt] {$f_2$};
		
					\draw (-.5, -.0) -- ( 1.5, -.0); 
					\fill (-.5, -.0) circle (2pt) node[below =3pt] {$-3$};
					\fill ( 1.5, -.0) circle (2pt) node[below=3pt] {$-2$};
					\fill (-.5, -.0) circle (2pt) node[above =3pt] {$f_3$};
					\fill ( 1.5, -.0) circle (2pt) node[above=3pt] {$f_4$};
				\end{tikzpicture}
			\]
			\caption{Plumbing 4-manifold, whose boundary is the 2-fold cover of $K\# K$, 
						with our chosen basis for its intersection form.}
			\label{figure eight plumbing} 
		}
		\end{figure}
	\end{example}

\section{Montesinos Knots and Their 2-Fold Covers}\label{mont def}

\begin{figure}[htp]
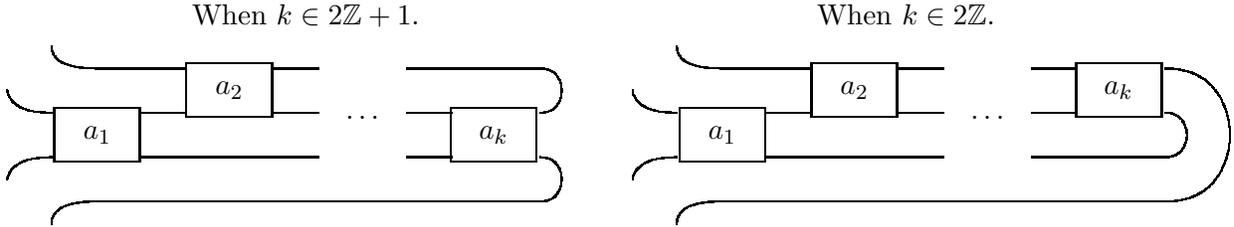

		\[
			\xygraph{	
				!{0;/r1.4pc/:}	
				!{\xbendr@(-.5)} !{\xcaph@(-.5)} !{\xcaph@(-.5)} [rr]
				!{\xcaph @(-.5)} [rr] !{\xcaph@(-.5)} !{\xcaph@(-.5)} 
				!{\xcaph@(-.5)} !{\hcap}  
				[llllllllllll]
				!{\xbendr@(-.4)} [rr] !{\xcaph@(-.5)} [rr]
				!{\xcaph@(-.5)} [rr]  !{\xcaph@(-.5)}
				[dllllllllll] 
				!{\xbendl-@(-.4)} [rrru] !{\xcaph@(-.5)} !{\xcaph@(-.5)} 
				!{\xcaph@(-.5)} !{\xcaph@(-.5)} [rr]  !{\xcaph@(-.5)}
				[rr]  !{\hcap}  
				[dlllllllllll]
				!{\xbendl-@(-.5)} [ur] !{\xcaph@(-.5)} !{\xcaph@(-.5)}
				!{\xcaph@(-.5)} !{\xcaph@(-.5)} !{\xcaph@(-.5)}
				!{\xcaph@(-.5)} !{\xcaph@(-.5)} !{\xcaph@(-.5)}
				!{\xcaph@(-.5)} !{\xcaph@(-.5)} 
				!{\put(18,-52) {\framebox(32,20){$a_1$}}}
				!{\put(68,-35) {\framebox(32,20){$a_2$}}}
				!{\put(128,-36) {$\ldots$}}
				!{\put(168,-52){\framebox(32,20){$a_k$}}}
				!{\put(70,0) {When $k\in 2\Z+1$.}}
			}
	 		\hspace{.25in}
			\xygraph{	
				!{0;/r1.4pc/:}	
				!{\xbendr@(-.5)} !{\xcaph@(-.5)} !{\xcaph@(-.5)} [rr]
				!{\xcaph @(-.5)} [rr] !{\xcaph@(-.5)} [rr] !{\hcap[3]}  
				[llllllllllll]
				!{\xbendr@(-.4)} [rr] !{\xcaph@(-.5)} [rr]
				!{\xcaph@(-.5)} [rr]  !{\xcaph@(-.5)} [rr] !{\hcap}  
				[dllllllllllll] 
				!{\xbendl-@(-.4)} [rrru] !{\xcaph@(-.5)} !{\xcaph@(-.5)} 
				!{\xcaph@(-.5)} !{\xcaph@(-.5)} [rr]  
				!{\xcaph@(-.5)} !{\xcaph@(-.5)} !{\xcaph@(-.5)}
				[dlllllllllll]
				!{\xbendl-@(-.5)} [ur] !{\xcaph@(-.5)} !{\xcaph@(-.5)}
				!{\xcaph@(-.5)} !{\xcaph@(-.5)} !{\xcaph@(-.5)}
				!{\xcaph@(-.5)} !{\xcaph@(-.5)} !{\xcaph@(-.5)}
				!{\xcaph@(-.5)} !{\xcaph@(-.5)} 
				!{\put(18,-52) {\framebox(32,20){$a_1$}}}
				!{\put(68,-35) {\framebox(32,20){$a_2$}}}
				!{\put(128,-36) {$\ldots$}}
				!{\put(168,-35){\framebox(32,20){$a_k$}}}
				!{\put(70,0) {When $k\in 2\Z$.}}
			}
		\]
		\caption{Rational tangles corresponding to
				 $\a/\b=[a_1,a_2,\ldots, a_k]$ for $k$ 
				 odd and even respectively.}
		\label{rational tangle diagram}
	\end{figure}

	Here, we outline the definition of a Montesinos knot and its 2-fold branched cover.  
	We follow \cite{Owens}.  Let $n_1, n_2,\ldots, n_k\in \Z$, then we let 
	$[n_1,n_2,\ldots, n_k]$ be the continued fraction
		\[
			[n_1,n_2,\ldots, n_k] = 
			n_1 - \frac{1}{\displaystyle n_2 - \frac{1}{\displaystyle \ldots - \frac{1}{n_k}}}.
		\]
	Let $\a/\b\in \Q$ be in lowest terms and let $[a_1,a_2,\ldots,a_k]$ be a continued 
	fraction representing $\a/\b$, then a rational tangle corresponding to $\a/\b$ is given 
	by Figure \ref{rational tangle diagram} where \framebox{$a_i$} indicates $a_i$ half-twists.  
	With these conventions in place, we can state the definition of a Montesinos link
	and describe its 2-fold cover.
	\begin{definition}\label{mont knots}
		Let $\a_1/\b_1, \a_2/\b_2,\ldots, \a_n/\b_n\in \Q$ each be in lowest terms and
		let $e\in\Z$.  Then the Montesinos link $M(e;\a_1/\b_1,\ldots, \a_n/\b_n)$ is
		given by Figure \ref{Montesinos link} where \framebox{$e$} represents $e$ 
		half-twists and \framebox{\tiny{$\a_i/\b_i$}} represent a rational 
		tangle corresponding to $\a_i/\b_i$.  
	\end{definition}
	\begin{figure}[htp]
		\[
			\xygraph
			{	
				!{0;/r1.4pc/:}
				!{\hcap[-4]} [d] !{\hcap[-2]} 
				[rr] !{\xcaph@(-.5)} [rr] !{\xcaph@(-.5)} [rrrr] !{\xcaph@(-.5)} [rr] !{\hcap[2]}
				[lllllllllllllu]
				[rr] !{\xcaph-@(-.5)} [rr] !{\xcaph-@(-.5)} [rrrr] !{\xcaph-@(-.5)} [rr] !{\hcap[4]}
				[lllllllllllllddd]
				!{\xcaph[-1]@(-.5)} !{\xcaph[-1]@(-.5)} !{\xcaph[-1]@(-.5)} 
				!{\xcaph[-1]@(-.5)} !{\xcaph[-1]@(-.5)} !{\xcaph[-1]@(-.5)}
				[rr] 
				!{\xcaph[-1]@(-.5)} !{\xcaph[-1]@(-.5)} !{\xcaph[-1]@(-.5)}
				!{\xcaph[-1]@(-.5)} !{\xcaph[-1]@(-.5)}
				[llllllllllllld]
				!{\xcaph[1]@(-.5)} !{\xcaph[1]@(-.5)} !{\xcaph[1]@(-.5)}
				!{\xcaph[1]@(-.5)} !{\xcaph[1]@(-.5)} !{\xcaph[1]@(-.5)}
				[rr] 
				!{\xcaph[1]@(-.5)} !{\xcaph[1]@(-.5)} !{\xcaph[1]@(-.5)}
				!{\xcaph[1]@(-.5)} !{\xcaph[1]@(-.5)}				!{\put(18,-18) {\framebox(32,20){$\a_1/\b_1$}}}
				!{\put(68,-18) {\framebox(32,20){$\a_2/\b_2$}}}
				!{\put(145,-10){\large{$\mathbf{\ldots}$}}}
				!{\put(202,-18){\framebox(32,20){$\a_n/\b_n$}}}
				!{\put(118,-69){\framebox(32,20){$e$}}}
			}
		\]
		\caption{Montesinos link $M(e;\a_1/\b_1,\ldots, \a_n/\b_n)$.}
		\label{Montesinos link}
	\end{figure}
	\begin{thm}[\cite{Montesinos}]\label{2-fold cover of plumbing proposition}
		The 2-fold branched cover of $S^3$ branched along a Montesinos link 
		$M(e;\a_1/\b_1,\ldots, \a_n/\b_n)$ is the boundary of the plumbing 4-manifold 
		given by Figure \ref{general plumbing}
		\begin{figure}[htp]
		{
			\[
				\begin{tikzpicture}
				{
					\draw[rotate=90] (0,0) -- (0,-2) (0,-3) -- (0,-3.25);
					\draw[rotate=100] (0,0) -- (2,0) (3,0) -- (3.25,0);
					\draw[rotate=180] (0,0) -- (2,0) (3,0) -- (3.25,0);
					\draw[rotate=140] (0,0) -- (2,0) (3,0) -- (3.25,0);
					\fill[rotate=180] 	(1,0) circle (2pt) node[above=0pt] {$a_1^1$}
										(1.75,0) circle (2pt) node[above=3pt] {$a_2^1$}
										(3.25  ,0) circle (2pt) node[above=3pt] {$a_{\ell_1}^1$}
										(2.25,0) circle (1pt)
										(2.50,0) circle (1pt)
										(2.750,0) circle (1pt);
					\fill[rotate=140] 	(1,0) circle (2pt) node[above=3pt] {$a_1^2$}
										(1.75,0) circle (2pt) node[above=3pt] {$a_2^2$}
										(3.25  ,0) circle (2pt) node[above=3pt] {$a_{\ell_2}^2$}
										(2.25,0) circle (1pt)
										(2.50,0) circle (1pt)
										(2.750,0) circle (1pt);	
					\fill[rotate=90]	(0,-1) circle (2pt) node[below=3pt] {$a_1^n$}
										(0,-1.75) circle (2pt) node[below=3pt] {$a_2^n$}
										(0,-3.25)   circle (2pt) node[below=3pt] {$a_{\ell_n}^n$}
										(0,-2.25) circle (1pt)
										(0,-2.50) circle (1pt)
										(0,-2.75) circle (1pt);
					\fill[rotate=100] 	(1,0) circle (2pt) node[right=3pt] {$a_1^3$}
										(1.75,0) circle (2pt) node[right=3pt] {$a_2^3$}
										(3.25  ,0) circle (2pt) node[right=3pt] {$a_{\ell_3}^3$}
										(2.25,0) circle (1pt)
										(2.50,0) circle (1pt)
										(2.750,0) circle (1pt);
					\fill (0,0) circle (2pt) node[below=3pt] {$e$};
					\fill[rotate=60]		(1.75,0) circle (1pt);
					\fill[rotate=45]		(1.75,0) circle (1pt);
					\fill[rotate=30]		(1.75,0) circle (1pt);
				}
				\end{tikzpicture}
			\]
			\caption{Plumbing 4-manifold $X_K$ with $\partial X_K = Y_K$ for $K=M(e;\a_1/\b_1,\ldots, \a_n/\b_n)$.}
			\label{general plumbing}
		}
		\end{figure}
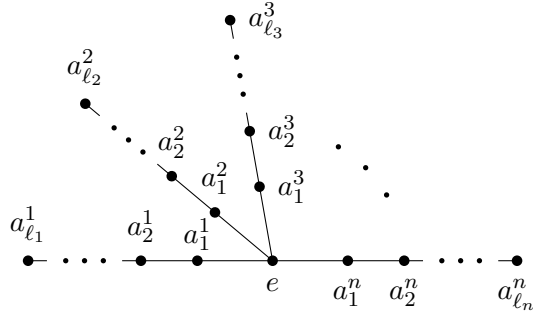
		where $\a_i/\b_i = [a_1^i,a_2^i,\ldots,a_{\ell_i}^i]$.
	\end{thm}
	Neumann and Raymond \cite{Neumann-Raymond} give the following useful result to
	determine when such a plumbing could be blown down to a negative definite 4-manifold.
	\begin{thm}[\cite{Neumann-Raymond}]\label{negative definiteness condition}
		A plumbing, $W$, with no nonnegative weights in its plumbing graph and whose 
		boundary is the 2-fold branched cover of $S^3$ branched over
		the Montesinos link $M(e;\a_1/\b_1,\ldots, \a_n/\b_n)$, is negative definite if
		and only if
		\[
			e + \sum_{i=1}^n \frac{\b_i}{\a_i} > 0.
		\]
	\end{thm}

	\section{Obstructing Sliceness in $M(0;[m_1+1,n_1+2],[m_2+1,n_2+2],q)$}\label{main prof section}
	
	In this section, we consider $M(0;[m_1+1,n_1+2],[m_2+1,n_2+2],q)$ 
	and prove Theorem \ref{main}.  First, we outline the restrictions on $m_i, n_i, q$ found
	within Theorem \ref{main}.  In order for  $M(0;[m_1+1,n_1+2],[m_2+1,n_2+2],q)$ to be a knot
	we need to carefully restrict the choices of parity on $m_i, n_i, q$.  The following
	proposition outlines these restrictions.
	\begin{proposition}\label{m_i n_i q parity proposition}
		All but the following twelve combinations of parity for the $m_i$, $n_i$, and $q$,
		result in knots in the family $M(0;[m_1+1,n_1+2],[m_2+1,n_2+2],q)$:
		\begin{enumerate}
			\item $q$, $m_1$, and $n_1+1$ even and no other restrictions,
			\item $q$, $m_2$, and $n_2+1$ even and no other restrictions,
			\item $q$, $m_1$, $n_1$, $n_2+1$ odd and no other restrictions,
			\item $q$, $n_1+1$, $m_2$, $n_2$ odd and no other restrictions,
			\item $q$, $m_1+1$, $n_1$, $m_2+1$, and $n_2$ odd.
		\end{enumerate}
	\end{proposition}
	\begin{proof}
		The proof is a simple matter of an exhaustive check of the thirty-two possible
		combinations of parity for the $m_i$, $n_i$ and $q$.  
	\end{proof}
	To apply Theorem \ref{main obstruction}, we need the 2-fold cover of $M(0;[m_1+1,n_1+2],[m_2+1,n_2+2],q)$
	to bound a negative definite 4-manifold.  In light of Theorem \ref{negative definiteness condition},
	we require that $m_i$, $n_i$ and $q$ satisfy
	\begin{equation}\label{m_i n_i restriction}
		\frac{n_1+2}{\left(m_1+1\right)n_1+2m_1+1}+\frac{n_2+2}{\left(m_2+1\right)n_2+2m_2+1}+\frac{1}{q}>0.
	\end{equation}

	\subsection{Choosing a Plumbing 4-manifold}
\begin{figure}[htp]
	{
		\[
			\begin{tikzpicture} 
				\draw (-2.4, 0) -- ( 2.4, 0); 
				\draw (-5.5, 0) -- (-3.35, 0);
				\draw (-7   , 0) -- (-6.75, 0); 
				\draw ( 5.5, 0) -- (3.35, 0);
				\draw ( 6.75, 0) -- ( 7   , 0);
				\draw ( 0   , 0) -- ( 0   ,-1);
				
				\draw[snake=brace] (-5.1,-.25) -- (-7.15,-.25) node[midway,below=3pt]{$n_1$};
				\draw[snake=brace] (-2.1,-.25) -- (-3.65,-.25) node[midway,below=3pt]{$m_1$};
				\draw[snake=brace] ( 7.15,-.25) -- ( 5.1,-.25) node[midway,below=3pt]{$n_2$};
				\draw[snake=brace] ( 3.65,-.25) -- ( 2.1,-.25) node[midway,below=3pt]{$m_2$};
				
				\fill (-6.425,0) circle (1pt);
				\fill (-6.125,0) circle (1pt);
				\fill (-5.825,0) circle (1pt); 
				
				\fill (-2.575,0) circle (1pt);
				\fill (-2.875,0) circle (1pt);
				\fill (-3.175,0) circle (1pt); 
				
				\fill ( 6.425,0) circle (1pt);
				\fill ( 6.125,0) circle (1pt);
				\fill ( 5.825,0) circle (1pt); 
				
				\fill ( 2.575,0) circle (1pt);
				\fill ( 2.875,0) circle (1pt);
				\fill ( 3.175,0) circle (1pt); 
				
	
				\fill ( 0, 0) circle (2pt) node[above=3pt] {0};
				\fill ( .75, 0) circle (2pt) node[above=3pt] {-1};
				\fill (-.75, 0) circle (2pt) node[above=3pt] {-1};
				\fill (-1.5,0) circle (2pt) node[above=3pt] {-1};
				\fill (1.5,0) circle (2pt) node[above=3pt] {-1};
				\fill (-2.25, 0) circle (2pt) node[above=3pt] {-2};
				\fill (2.25, 0) circle (2pt) node[above=3pt] {-2};
		
				\fill (-5.25, 0) circle (2pt) node[above=3pt] {-2};
				\fill (-7, 0) circle (2pt) node[above=3pt] {-2};
				\fill ( 0,-1) circle (2pt) node[below=3pt] {$q$};
				\fill ( 3.5, 0) circle (2pt) node[above=3pt] {-2};
				\fill ( -3.5, 0) circle (2pt) node[above=3pt] {-2};
				\fill ( 4.375, 0) circle (2pt) node[above=3pt] {-3};
				\fill ( -4.375, 0) circle (2pt) node[above=3pt] {-3};
				\fill ( 5.25, 0) circle (2pt) node[above=3pt] {-2};
				\fill ( 7, 0) circle (2pt) node[above=3pt] {-2};
			\end{tikzpicture}
		\]
		\caption{Plumbing 4-manifold whose boundary is $Y_K$ for $K=M(0;[m_1+1,n_1+2],[m_2+1,n_2+2],q)$.}
		\label{plumbing graph pre bd}
	}
	\end{figure}
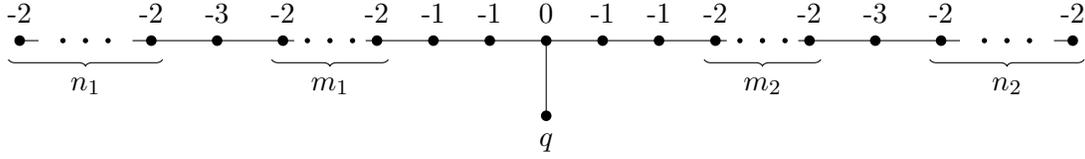

	Rather than applying Theorem \ref{2-fold cover of plumbing proposition} 
	directly, we note that $[m+1,n+2]$ can be expressed as a different continued fraction expansion
	 - one which is more compatible with Lisca's approach.
	\begin{proposition}\label{CF result}
		Let $m, n\in \N$, then  
		$[m+1,n+2]=\displaystyle [-1,-1,\underbrace{-2,\ldots, -2}_{m},-3,\underbrace{-2,\ldots, -2}_{n}]$.
	\end{proposition}
	\begin{proof}
		This fact is easily proven by induction.
	\end{proof}
	This, combined with Theorem \ref{2-fold cover of plumbing proposition}, gives that 
	the 2-fold cover of $M(0;[m_1+1,n_1+2],[m_2+1,n_2+2],q)$ bounds the plumbing 4-manifold
	in Figure \ref{plumbing graph pre bd}.  Blowing down appropriate vertices results in the
	negative definite 4-manifold in Figure \ref{plumbing graph}.
	
	\subsection{Proving Theorem \ref{main}}\label{main proof}
	
		Before we present the proof of Theorem \ref{main}, we define our notation for
		embeddings of intersection forms over the integers.  
		For an integer $n$, we view the $n\times n$ matrix $A$ as an intersection pairing 
		on $\Z^n$ with respect to the standard basis $f_1,\ldots, f_n$.  Similarly, we view $-I$ as
		the standard negative definite intersection pairing on another copy of $\Z^n$
		with respect to its standard basis, $e_1,\ldots, e_n$.  We abbreviate $\<Af_i,f_j\>$
		to $f_i \cdot f_j$, and similarly for $-I$.   Finally, an embedding of $A$ into 
		$-I$ is a monomorphism $\varphi:\Z^n\to \Z^n$ such that $f_i\cdot f_j=
		\varphi(f_i)\cdot \varphi (f_j)$.   
		

	\begin{proof}[Proof of Theorem \ref{main}]
		Let $Q_W$ be the intersection form associated to the plumbing 4-manifold in
		Figure \ref{plumbing graph}, we show that $Q_W$ does not embed into the standard negative 
		definite form of equal rank.
		\begin{figure}[htp]
		{
			\[
			\begin{tikzpicture}[scale=.95]
				\draw (-1.25, 0) -- ( 1.25, 0); 
				\draw (-5.25, 0) -- (-2.75, 0);
				\draw (-7   , 0) -- (-6.75, 0); 
				\draw ( 2.75, 0) -- ( 5.25, 0);
				\draw ( 6.75, 0) -- ( 7   , 0);
				\draw ( 0   , 0) -- ( 0   ,-2);
		
				\draw[snake=brace] (-4.85,-.75) -- (-7.15,-.75) node[midway,below=3pt]{$n_1$};
				\draw[snake=brace] (-0.85,-.75) -- (-3.15,-.75) node[midway,below=3pt]{$m_1-1$};
				\draw[snake=brace] ( 7.15,-.75) -- ( 4.85,-.75) node[midway,below=3pt]{$n_2$};
				\draw[snake=brace] ( 3.15,-.75) -- ( 0.85,-.75) node[midway,below=3pt]{$m_2-1$};
				
				\fill (-6.35,0) circle (1pt);
				\fill (-6.00,0) circle (1pt);
				\fill (-5.65,0) circle (1pt); 
				
				\fill (-2.35,0) circle (1pt);
				\fill (-2.00,0) circle (1pt);
				\fill (-1.65,0) circle (1pt); 
				
				\fill ( 6.35,0) circle (1pt);
				\fill ( 6.00,0) circle (1pt);
				\fill ( 5.65,0) circle (1pt); 
				
				\fill ( 2.35,0) circle (1pt);
				\fill ( 2.00,0) circle (1pt);
				\fill ( 1.65,0) circle (1pt); 
				
				\fill (-7, 0) circle (2pt) node[below=3pt] {-$2$} node[right=4pt,rotate=70,xshift=2pt,yshift=10pt] {$f_1$};
				\fill (-5, 0) circle (2pt) node[below=3pt] {-$2$} node[right=4pt,rotate=70,xshift=2pt,yshift=10pt] {$f_{n_1}$};
				\fill (-4, 0) circle (2pt) node[below=3pt] {-$3$} node[right=4pt,rotate=70,xshift=2pt,yshift=10pt] {$f_{n_1+1}$};
				\fill (-3, 0) circle (2pt) node[below=3pt] {-$2$} node[right=4pt,rotate=70,xshift=2pt,yshift=10pt] {$f_{n_1+2}$};
				\fill (-1, 0) circle (2pt) node[below=3pt] {-$2$} node[right=4pt,rotate=70,xshift=2pt,yshift=10pt] {$f_{n_1+m_1}$};
				\fill ( 0, 0) circle (2pt) node[below=3pt] {\hspace{.5cm} -$2$}  
								   node[right=4pt,rotate=70,xshift=2pt,yshift=10pt] {$f_{n_1+m_1+1}$};
				\fill ( 1, 0) circle (2pt) node[below=3pt] {-$2$} node[right=4pt,rotate=70,xshift=2pt,yshift=10pt] {$f_{n_1+m_1+2}$};
				\fill ( 3, 0) circle (2pt) node[below=3pt] {-$2$} node[right=4pt,rotate=70,xshift=2pt,yshift=10pt] {$f_{n_1+m_1+m_2}$};
				\fill ( 4, 0) circle (2pt) node[below=3pt] {-$3$} node[right=4pt,rotate=70,xshift=2pt,yshift=10pt] {$f_{n_1+m_1+m_2+1}$};
				\fill ( 5, 0) circle (2pt) node[below=3pt] {-$2$} node[right=4pt,rotate=70,xshift=2pt,yshift=10pt] {$f_{n_1+m_1+m_2+2}$};
				\fill ( 7, 0) circle (2pt) node[below=3pt] {-$2$} node[right=4pt,rotate=70,xshift=2pt,yshift=10pt] {$f_{n_1+m_1+m_2+n_2+1}$};
				\fill ( 0,-2) circle (2pt) node[left=3pt] {$q$} node[below=3pt] {$f_{n_1+m_1+m_2+n_2+2}$};
			\end{tikzpicture}			
			\]
			\caption{Plumbing 4-manifold with our chosen basis for $Q_W$.}
			\label{plumbing graph} 
		}
		\end{figure}
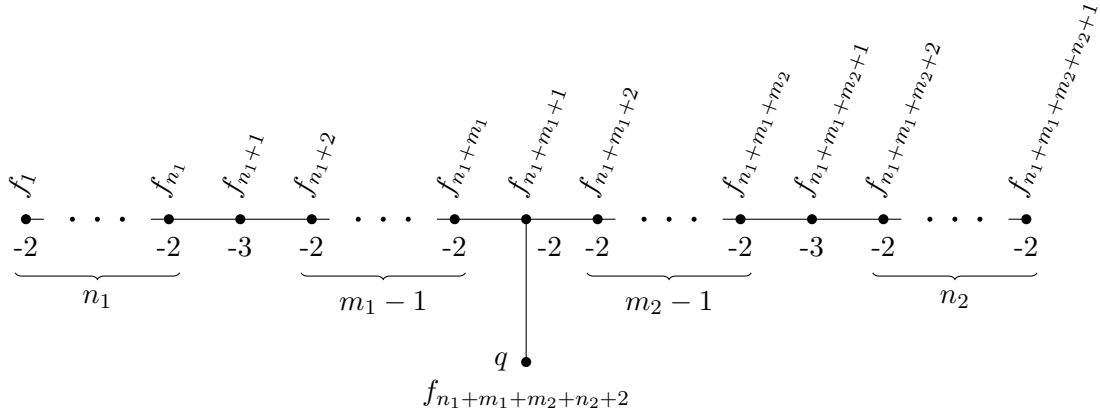
		Suppose to the contrary that $Q_W$ does embed into $-I$.  Then, there must exist 
		an embedding $\varphi:\Z^n\to \Z^n$, where $n=n_1+m_1+m_2+n_2+2$.  
		Taking the $f_i$ indicated in Figure \ref{plumbing graph}
		as our basis for $Q_W=[q_{ij}]$, we explore the structure of $\varphi$.  
		
		First, note that up to a change of basis for $-I$, any vertex with self intersection -2 or 
		-3 has essentially a unique image under $\varphi$.  For instance, we can assume that, 
		$\varphi\left(f_1\right) = e_1-e_2$. Similarly, we can assume that 
		$\varphi\left(f_2\right)=e_i-e_{i+1}$ for some $i\in\left\{1,2,\ldots,n\right\}$ (up to a change of basis).  
		Now $f_1\cdot f_2 = q_{12}=1$; therefore $\varphi\left(f_1\right)\cdot \varphi\left(f_2\right)=1$.  
		By definition		
		\[
			\varphi \left(f_1\right)\cdot \varphi\left(f_2\right) = e_1\cdot e_i-e_1\cdot e_{i+1}-e_2\cdot e_i +e_2\cdot e_{i+1}.
		\]
		Since $\varphi$ is injective, we cannot have $i=1$.  If $i>2$, then the above quantity 
		is zero, so it follows that $i=2$ and thus $\varphi\left(f_2\right)=e_2-e_3$.  Similar arguments 
		show that 
		\begin{equation}\label{pairing images}
			\varphi\left(f_i\right) = e_i-e_{i+1}, \hspace{.25in} 1\leq i\leq n_1.
		\end{equation}
		It is worth noting that if $n_1>2$, we could take $\varphi(f_3)=-e_1-e_2$ while still 
		satisfying the pairings $f_1\cdot f_3=0$ and $f_2\cdot f_3=1$.  However, with such 
		a definition, we still need $f_4\cdot f_3=1$, implying that $\varphi(f_4)$ must contain either 
		$\pm e_1$ or $\pm e_2$, but not both.  However, this would imply that $\varphi(f_1)$ would 
		pair with $\varphi(f_4)$ nontrivially - which contradicts the fact that $\varphi$ is an embedding.
		It follows that, up to a change of basis for $-I$, the images
		defined in \eqref{pairing images} are the only possibilities.  
		
		Next, we consider the image of $f_{n_1+2}$ under $\varphi$.  As above,  
		up to a change of basis we can take 
		$\varphi\left(f_{n_1+2}\right) = e_{n_1+2}-e_{n_1+3}$ and by the same argument used for 
		$f_2$, we arrive at
		\[
			\varphi\left(f_i\right) = e_i-e_{i+1}, \hspace{.25in} n_1+2\leq i\leq n_1+m_1+m_2.
		\]
		Similarly, we can assume that $\varphi\left(f_{n_1+m_1+m_2+2}\right) = e_{n_1+m_1+m_2+2}-e_{n_1+m_1+m_2+3}$ 
		and thus
		\[
			\varphi\left(f_i\right) = e_i-e_{i+1}, \hspace{.25in} n_1+m_1+m_2+2\leq i\leq n_1+m_1+m_2+n_2+1.
		\]

		Next, we consider $\varphi\left(f_{n_1+1}\right)$.  Since $\varphi$ is a homomorphism, we know that 
		$\varphi\left(f_{n_1+1}\right) = \sum_{i=1}^{n_1+m_1+m_2+n_2+2}\la_i e_i$ for some $\la_i\in \Z$.  
		Given that $f_i\cdot f_{n_1+1} = \varphi\left(f_i\right)\cdot \varphi\left(f_{n_1+1}\right)$, we have
		\begin{equation}\label{n_1+1 pairings}
			\varphi\left(f_i\right)\cdot \varphi\left(f_{n_1+1}\right) = \left\{ 	
				\begin{array}{ll}
					\pm 1	&  	i= n_1,		\\
					-3 		&	i=n_1+1,	\\
					1		& 	i=n_1+2,	\\
					0		&	\mbox{otherwise.}
				\end{array}
			\right.
		\end{equation}
		\noindent Note that the $\pm1$ (rather than just 1) arises from the fact that we may have had to 
		change the basis of $-I$ to get images of the previous vertices in their correct forms.  
		Thus, we may have caused the pairing of $\varphi\left(f_{n_1}\right)$ and $\varphi\left(f_{n_1+1}\right)$ 
		to become negative.  
		
		Since $\varphi\left(f_i\right) = e_i-e_{i+1}$ for each $i \neq n_1+1, n_1+m_1+m_2+1, 
		n_1+m_1+m_2+n_2+2$, we have that  $\varphi\left(f_i\right)\cdot \varphi\left(f_{n_1+1}\right) 
		= \la_{i+1}-\la_i$ for the same $i$.  Then, from \eqref{n_1+1 pairings} we have
		that
		\[
			\la_{i+1}-\la_i = 0,
		\]
		for $1\leq i\leq n_1-1$, $n_1+3\leq i \leq n_1+m_1+m_2$, and $n_1+m_1+m_2+3\leq n_1+m_1+m_2+n_2+1$.  
		It follows that 
		\[
			\la_1=\la_2=\ldots = \la_{n_1}, \hspace{.25in}	\la_{n_1+3} = \la_{n_1+4} = \ldots = \la_{n_1+m_1+m_2+1},
		\]
		\[
			\la_{n_1+m_1+m_2+3} = \la_{n_1+m_1+m_2+4} = \ldots = \la_{n_1+m_1+m_2+n_2+2}.
		\]
		Let $\la = \la_1$, $\eta = \la_{n_1+3}$ and $\g = \la_{n_1+m_1+m_2+3}$.  Now, 
		$\varphi\left(f_{n_1}\right)\cdot \varphi\left(f_{n_1+1}\right) = \la_{n_1}-\la_{n_1+1} = \la_{n_1}-\la = \pm 1$. 
		Therefore, $\la_{n_1} = \la \pm 1$.  Similarly, $\la_{n_1+2} = \eta - 1$.  So,
		\begin{align}\nonumber
			\varphi\left(f_{n_1+1}\right) &=\cr
			&\hspace{-15mm}=\sum_{i=1}^{n_1}\la e_i+\left(\la \pm 1\right)e_{n_1+1}+\left(\eta -1\right)e_{n_1+2} 	
			+\sum_{i=n_1+3}^{n_1+m_1+m_2+1}\eta e_i+\sum_{i=n_1+m_1+m_2+2}^{n_1+m_1+m_2+n_2+2}\g e_i.
		\end{align}
		Direct calculation along with \eqref{n_1+1 pairings} gives that
		\[
			\varphi\left(f_{n_1+1}\right)^2 = 
				-\left(n_1 \la^2 + \left(\la\pm 1\right)^2+\left(\eta - 1\right)^2 +\left(m_1+m_2-1\right)\eta^2+
					\left(n_2+1\right)\g^2\right)=-3.
		\]
		Clearly, if $n_2>2$, then $\g=0$.  If $n_2=2$ and $\g\neq 0$, then $\la^2$, 
		$\left(\la\pm 1\right)^2$, $\left(\eta - 1\right)^2$, and $\eta^2$ must be be identically zero - 
		clearly an impossibility; therefore if $n_2=2$, $\g=0$. Moreover, if $n_2=1$ and 
		$\g\neq 0$, then three of $\la^2$, $\left(\la\pm 1\right)^2$, $\left(\eta - 1\right)^2$, and 
		$\eta^2$ must be zero - which is again impossible.  It follows that $\g$ is necessarily 0.  
		Therefore
		\[
			\varphi\left(f_{n_1+1}\right) = 
				\sum_{i=1}^{n_1}\la e_i+\left(\la \pm 1\right)e_{n_1+1}+\left(\eta - 1\right)e_{n_1+2}+
				\sum_{i=n_1+3}^{n_1+m_1+m_2+1}\eta e_i.
		\]
		The same argument applies to $\varphi\left(f_{n_1+m_2+m_2+1}\right)$ where the pairings are given 
		by
		\begin{equation}\label{n_1+m_1+m_2+1 pairings}
			\varphi\left(f_i\right)\cdot \varphi\left(f_{n_1+m_1+m_2+1}\right) =	\left\{ 	
				\begin{array}{ll}
					\pm 1	&  i= n_1+m_1+m_2+2,\\
					-3 		&	 i=n_1+m_1+m_2+1,\\
					1		& i=n_1+m_1+m_2,\\						0		&	\mbox{otherwise.}
				\end{array}
			\right.
		\end{equation}
		
		\noindent Then, there exist $\a,\b\in \Z$ such that 
		\begin{align}\nonumber
			\varphi \left( f_{n_1+m_1+m_2+1}\right) &= \cr
			& \hspace{-28mm}= \sum_{i=n_1+2}^{n_1+m_1+m_2}\b e_i+\left(\b + 1\right)e_{n_1+m_1+m_2+1}+ \left(\a\pm 1\right)e_{n_1+m_1+m_2+2} 
					+\sum_{i=n_1+m_1+m_2+3}^{n_1+m_1+m_2+n_2+2}\a e_i.
		\end{align}
		Therefore, we have explicit forms for the images of all but the last basis element in $Q_W$
		under $\varphi$.  Before we consider this last image, we extract as much information 
		out of the pairings in \eqref{n_1+1 pairings} and \eqref{n_1+m_1+m_2+1 pairings} as we can.  
		\eqref{n_1+m_1+m_2+1 pairings} gives that $\varphi\left(f_{n_1+1}\right)\cdot 
		\varphi\left(f_{n_1+m_1+m_2+1}\right) = 0$.  Therefore,
		\[
			\varphi\left(f_{n_1+1}\right)\cdot \varphi\left(f_{n_1+m_1+m_2+1}\right) 
				= -\left(\b\left(\eta - 1\right)+\left(m_1+m_2-2\right)\b\eta +\left(\b + 1\right)\eta\right)
				=0
		\]
		and thus
		\[
			\left(m_1+m_2\right)\b\eta - \b + \eta = 0.
		\]
		Noting that the coefficient of each 
		basis element in the image of a square -3 vertex is necessarily in the set $\{0,-1,1\}$, 
		we have that $\b, \b+1, \eta, \eta-1 \in \left\{0,-1,1\right\}$; so, $\b=0$ or $\beta=-1$ 
		and $\eta=0$ or $\eta=1$.

		If $\beta=0$, then $\eta=0$ and if $\beta=-1$, then $\eta=1$ and $m_1+m_2=2$.  
		In this latter case, $\a$ is forced to be $\pm 1$ which, in turn, forces $n_2=2$.  
		Moreover, $\la$ is forced to be $\pm1$ which makes $n_1=2$.  We'll come back to 
		this case.  
		
		Now, suppose that $\b = \eta =0$, then
		\[
			\varphi\left(f_{n_1+1}\right) = \sum_{i=1}^{n_1}\la e_i+\left(\la\pm 1\right)e_{n_1+1} - e_{n_1+2},
		\]
		\[
			\varphi\left(f_{n_1+m_1+m_2+1}\right)= 
				e_{n_1+m_1+m_2+1}+ \left(\a\pm 1\right)e_{n_1+m_1+m_2+2}+\sum_{i=n_1+m_1+m_2+3}^{n_1+m_1+m_2+n_2+2}\a e_i.
		\]	 
		Moreover, $\la = \mp 1$ and $n_1=2$, $\a=\mp 1$ and $n_2=2$, so
		\[
			\varphi\left(f_{n_1+1}\right) = \varphi\left(f_{3}\right) = \mp e_1 \mp e_2 -e_{4},
		\]
		\[
			\varphi\left(f_{n_1+m_1+m_2+1}\right) = \varphi\left(f_{m_1+m_2+3}\right)= e_{m_1+m_2+3} \mp e_{+m_1+m_2+5} \mp e_{m_1+m_2+6}.
		\]	 
		Therefore, if $n_1$ and $n_2$ are not identically 2, then no such embedding $\varphi$ exists.  
		
		Corresponding to the case when $\beta=\eta=0$, we have $n_1=n_2=2$ and $m_1, m_2\in \N$.
		Then, we know that up to a change of basis, 
		\begin{eqnarray*}
			\varphi\left(f_1\right) &=& e_1-e_2\\ 
			\varphi\left(f_2\right) &=& e_2-e_3\\ 
			\varphi\left(f_3\right) &=& \mp e_1 \mp e_2 -e_{4}\\ 
			\varphi\left(f_4\right) &=& e_4-e_5\\
							&\vdots&\\
			\varphi\left(f_{m_1+m_2+2}\right) &=& e_{m_1+m_2+3}-e_{m_1+m_2+2}\\ 
			\varphi\left(f_{m_1+m_2+3}\right) &=& e_{m_1+m_2+3} \mp e_{+m_1+m_2+5} \mp e_{m_1+m_2+6}\\ 
			\varphi\left(f_{m_1+m_2+4}\right) &=& e_{m_1+m_2+5}-e_{m_1+m_2+4}\\ 
			\varphi\left(f_{m_1+m_2+5}\right) &=& e_{m_1+m_2+6}-e_{m_1+m_2+5}
		\end{eqnarray*}
		Now, suppose that $\varphi\left(f_{m_1+m_2+6}\right) = \sum_{i=1}^{m_1+m_2+6} \mu_i e_i$ for $\mu_i\in \Z$.  
		Then, we get the following system of equations arising from the pairings indicated in
		Figure \ref{plumbing graph}
		\[
			\left\{ 	
			\begin{array}{l}	
				\varphi\left(f_1\right)\cdot \varphi\left(f_{m_1+m_2+6}\right) = \mu_2-\mu_{1} = 0 \\
				\varphi\left(f_2\right)\cdot \varphi\left(f_{m_1+m_2+6}\right) = \mu_3-\mu_{2} = 0 \\
				\varphi\left(f_i\right)\cdot \varphi\left(f_{m_1+m_2+6}\right) = \mu_{i+1}-\mu_{i} = 0 
						\hspace{.2in} 4 \leq i \leq m_1+2, \hspace{.15in}  m_1+4\leq i\leq m_1+m_2+2\\ 
				\varphi\left(f_{m_1+m_2+4}\right)\cdot \varphi\left(f_{m_1+m_2+6}\right) = \mu_{m_1+m_2+5}-\mu_{m_1+m_2+5} = 0  \\
				\varphi\left(f_{m_1+m_2+5}\right)\cdot \varphi\left(f_{m_1+m_2+6}\right) = \mu_{m_1+m_2+6}-\mu_{m_1+m_2+5} = 0  \\
				\varphi\left(f_{m_1+3}\right)\cdot \varphi\left(f_{m_1+m_2+6}\right) = \mu_{m_1+4}-\mu_{m_1+3} = 1  \\
				\varphi\left(f_{3}\right)\cdot \varphi\left(f_{m_1+m_2+6}\right) = \pm\mu_1\pm \mu_2 +\mu_4 = 0  \\
				\varphi\left(f_{m_1+m_2+3}\right)\cdot \varphi\left(f_{m_1+m_2+6}\right) 
						= -\mu_{m_1+m_2+3} \pm \mu_{+m_1+m_2+5} \pm \mu_{m_1+m_2+6} = 0 \\
			\end{array}
			\right.
		\]
		Then $\mu_1=\mu_2=\mu_3$, $\mu_4=\mu_5=\ldots = \mu_{m_1+3}$, $\mu_{m_1+4}=\mu_{m_1+5} = 
		\ldots = \mu_{m_1+m_2+3}$, and $\mu_{m_1+m_2+4}=\mu_{m_1+m_2+5}=\mu_{m_1+m_2+6}$. 
		Let $\mu=\mu_1$, $\nu=\mu_4$, $\rho = \mu_{m_1+4}$ and $\s = \mu_{m_1+m_2+4}$, then
		\begin{align}\nonumber
			\varphi\left(f_{m_1+m_2+6}\right) &=\cr 
				&\hspace{-15mm}= \mu\left(e_1+e_2+e_3\right) +\sum_{i=4}^{m_1+3}\nu e_i + 
					\sum_{i=m_1+4}^{m_1+m_2+3}\rho e_i+\s\left(e_{m_1+m_2+4}+\ldots+e_{m_1+m_2+6}\right)
		\end{align}
		By construction, the above satisfies the first five equations in the system.  
		Since $e_{m_1+3}-e_{m_1+4} = 1$, we have that $\nu=\rho+1$; 
		now, $\pm\mu_1\pm \mu_2 +\mu_4 = 0$ gives that $\mp 2\mu = \nu$ and thus 
		$\nu$ is necessarily even which means that $\rho$ is odd; however, 
		$-\mu_{m_1+m_2+3} \pm \mu_{+m_1+m_2+5} \pm \mu_{m_1+m_2+6} = 0$ gives that 
		$\rho = \pm 2\s$ - a contradiction.  Therefore no such function $\varphi$ could exist if $\b=\eta$.
		
		It follows that the only way for such $\varphi$ to have any hope of
		existing is if we fall into the latter case mentioned above.  In this case,
		$m_1=m_2=1$, $n_1=n_2=2$, $\beta=-1$, $\eta=1$, $\a=\pm 1$, and $\la=\pm1$.
		Then, again, up to a change of basis, we can take
		\begin{eqnarray*}
			\varphi\left(f_1\right) &=& e_1-e_2\\ 
			\varphi\left(f_2\right) &=& e_2-e_3\\ 
			\varphi\left(f_3\right) &=& \pm e_1\pm e_2 + e_5\\ 
			\varphi\left(f_4\right) &=& e_4-e_5\\ 
			\varphi\left(f_5\right) &=& -e_4 \pm e_7 \pm e_8\\
			\varphi\left(f_6\right) &=& e_6-e_7\\ 
			\varphi\left(f_7\right) &=& e_7-e_8
		\end{eqnarray*}
		Suppose that $\varphi\left(f_8\right) = \sum_{i=1}^8 \mu_i e_i$ for $\mu_i\in \Z$.  
		As above, we get the following system of equations arising from the 
		pairings indicated in Figure \ref{plumbing graph}
		\[
			\begin{array}{lcl}
				\varphi\left(f_1\right)\cdot \varphi\left(f_8\right) = \mu_2-\mu_1 = 0	
					&	\hspace{.25in}	
					&	\varphi\left(f_5\right)\cdot \varphi\left(f_8\right) = \mu_4\mp\mu_7\mp \mu_8 = 0\\
				\varphi\left(f_2\right)\cdot \varphi\left(f_8\right) = \mu_3-\mu_2 = 0	
					&	\hspace{.25in}	
					&	\varphi\left(f_6\right)\cdot \varphi\left(f_8\right) = \mu_7-\mu_6 = 0\\
				\varphi\left(f_3\right)\cdot \varphi\left(f_8\right) = \mp\mu_1\mp\mu_2-\mu_5 = 0
					&	\hspace{.25in}	
					&	\varphi\left(f_7\right)\cdot \varphi\left(f_8\right) = \mu_8-\mu_7 = 0\\
				\varphi\left(f_4\right)\cdot \varphi\left(f_8\right) = \mu_5-\mu_4 = 1	
					&	\hspace{.25in}	&	
			\end{array}
		\]
		Let $\mu=\mu_1$ and $\nu=\mu_6$, then $\mu=\mu_1=\mu_2=\mu_3$ and $\nu=\mu_6=\mu_7=\mu_8$.  
		So $\mu_4 = \pm 2\nu$ and $\mu_5 = \mp 2\mu$.  But, $\mu_5=1+\mu_4$, which means that 
		$\mp 2\mu = 1 \pm 2\nu$ - a contradiction.  So no such function $\varphi$ exists in this case
		either.
		
		Given that we have exhausted all possible cases, it is clear that no embedding 
		$\varphi:\Z^n\to \Z^n$ could exist for $Q_W$ arising
		from the plumbing indicated in Figure \ref{plumbing graph}.  It follows that
		$Q_W$ does not embed into the standard negative definite form $-I$ as claimed.
		Then, the result follows from Theorem \ref{main obstruction}.  Therefore, no member
		of $M([m_1+1,n_1+2],[m_2+1,n_2+2],q)$ with the outlined restrictions
		on $m_i, n_i, q$ can be smoothly slice.  
	\end{proof}

\section{The Infinite Subfamily}\label{infinite subfamily section}
			\begin{proof}[Proof of Theorem \ref{infinite family}]
			Consider the case with $m_1=n_2=1$, $m_2=2$, and $q=-3$.  Then,
			we show that $M(0;[2,n_1+2],[3,3],-3)$ contains our desired subfamily.  With these
			choices of $m_i$, $n_2$, and $q$, we can take any choice of $n_1$ to meet both the 
			requirement for negative definiteness and the requirements on parity.  Thus every 
			knot $M(0;[2,n_1+2],[3,3],-3)$ meets the restrictions outlined in Theorem \ref{main}.
			
			\begin{figure}[htp]
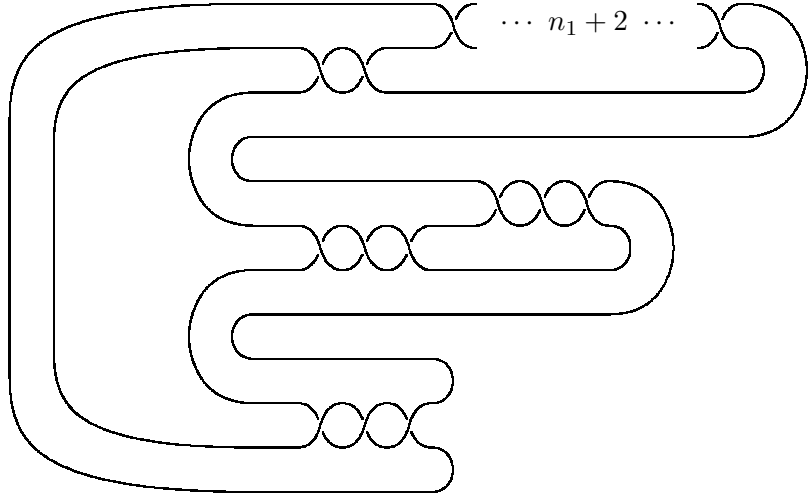

			{
			\[
			\xygraph
				{
				!{0;/r1.4pc/:}	
				!{\hcap[-11]@(.5)} !{\xcaph@(-.5)}
				!{\xcaph@(-.5)}	!{\xcaph@(-.5)} !{\xcaph@(-.5)}
				!{\htwist}[rrrrr] !{\htwist} !{\hcap[3]}
				[dlllllllllll]
				!{\xcaph@(-.5)} !{\htwist}!{\htwist}  !{\xcaph@(-.5)}
				[rrrrrrr] !{\hcap}
				[dlllllllllll]
				!{\xcaph@(-.5)}	[rr]
				!{\xcaph@(-.5)}	!{\xcaph@(-.5)} !{\xcaph@(-.5)}
				!{\xcaph@(-.5)}	!{\xcaph@(-.5)} !{\xcaph@(-.5)}
				!{\xcaph@(-.5)}	!{\xcaph@(-.5)} 
				[dlllllllllll]
				!{\xcaph@(-.5)}	!{\xcaph@(-.5)} !{\xcaph@(-.5)}
				!{\xcaph@(-.5)}	!{\xcaph@(-.5)} !{\xcaph@(-.5)}
				!{\xcaph@(-.5)}	!{\xcaph@(-.5)} !{\xcaph@(-.5)}
				!{\xcaph@(-.5)} !{\xcaph@(-.5)}
				[dlllllllllll]
				!{\xcaph@(-.5)}
				!{\xcaph@(-.5)}	!{\xcaph@(-.5)} !{\xcaph@(-.5)} !{\xcaph@(-.5)}  !{\htwist}!{\htwist} !{\htwist}
		 		!{\hcap[3]}
				[dllllllll]
				!{\xcaph@(-.5)} !{\htwist}!{\htwist} !{\htwist} !{\xcaph@(-.5)}
				[rrr] !{\hcap}
				[dllllllll]
				!{\xcaph@(-.5)}	[rrr]
				!{\xcaph@(-.5)}	!{\xcaph@(-.5)} !{\xcaph@(-.5)}!{\xcaph@(-.5)}
				[dllllllll]
				!{\xcaph@(-.5)}	!{\xcaph@(-.5)} !{\xcaph@(-.5)}
				!{\xcaph@(-.5)}	!{\xcaph@(-.5)} !{\xcaph@(-.5)}!{\xcaph@(-.5)} !{\xcaph@(-.5)}
				[uuuullllllll]
				!{\hcap-} [u] !{\hcap[-3]}
				[dddddd]
				!{\xcaph@(-.5)}	!{\xcaph@(-.5)}	!{\xcaph@(-.5)}!{\xcaph@(-.5)}
				!{\hcap}
				[dllll]
				!{\xcaph@(-.5)} 
				[dl]
				!{\xcaph@(-.5)}	[u] !{\htwistneg}!{\htwistneg} !{\htwistneg} [d]
				!{\hcap}
				[dllll]
				!{\xcaph@(-.5)} !{\xcaph@(-.5)}	!{\xcaph@(-.5)} !{\xcaph@(-.5)}
				[lllluuuu]
				!{\hcap-} [u] !{\hcap[-3]}
				[uuuuu]
				!{\hcap[-9]@(.5)}	
				!{\put(110,-10) {{$\cdots\ n_1+2\ \cdots$}}}
			}
			\]
				\caption{$M(0;[2,n_1+2],[3,3],-3)$}
				\label{family coloring}
			}
			\end{figure}
			
			Now, we use the method outlined in \cite{Gordon-Litherland} to calculate the 
			signature and determinant of each member
			of $M(0;[2,n_1+2],[3,3],-3)$.  We begin by checkerboard coloring the knot diagram 
			in Figure \ref{family coloring} with the unbounded exterior region colored white.  
			From this coloring, we arrive at the symmetric Goeritz matrix given by 
			\[
				G = 
				\begin{bmatrix}
					2  & -1 & 0 & 0 & 0 & 0 & 0 & 0\\
					-1 & -1 & 1 & 0 & 0 & 1 & 0 & 0\\
					0  & 1  &-2 & 1 & 0 & 0 & 0 & 0\\
					0  & 0  & 1 &-2 & 1 & 0 & 0 & 0\\
					0  & 0 & 0  & 1 & 2 & 0 & 0 &-3\\
					0  & 1 & 0  & 0 & 0 &-2 & 1 & 0\\
					0  & 0 & 0  & 0 & 0 & 1 & n_1+1 &-n_1-2\\
					0  & 0 & 0  & 0 &-3 & 0 &-n_1-2 & n_1+6
				\end{bmatrix}
			\]	
			Direct calculation gives that $\det(G) = 26n_1+51$, which as \cite{Gordon-Litherland} shows,
			is precisely the Determinant of the knot $M(0;[2,n_1+2],[3,3],-3)$.  Summing over all type II
			crossings in Figure \ref{family coloring} with our chosen coloring scheme and an arbitrary
			string orientation gives a correction factor of $\mu(M(0;[2,n_1+2],[3,3],-3)) = 3-3=0$.  Therefore,
			the signature of $M(0;[2,n_1+2],[3,3],-3)$ is precisely the signature of $G$.  It's easy to verify
			that the signature of $G$ is 0 for any valid choice of $n_1$ and thus
			each member of $M(0;[2,n_1+2],[3,3],-3)$ has signature 0.  
			
			Therefore, each $n_1$ such that $26n_1+51=a^2$ for some $a\in \N$, gives an
			example of a knot in the family $M(0;[m_1+1,n_1+2],[m_2+1,n_2+2],q)$ which
			has square determinant and signature 0.  To show that this is an infinite subset, 
			consider the sequence $\left\{a_n\right\}_{n=0}^\infty$
			given recursively by $a_0=21$, and
			\[
				a_{n+1} = 	\left\{\begin{array}{ll}
									a_n + 10 & \mbox{if $n$ is even,}\\
									a_n + 16 & \mbox{if $n$ is odd.}
								\end{array}
							\right.
			\]
			\noindent We claim that there is a positive integral solution, $n_1$, to $26n_1+51=a_n^2$
			for each $n$.  To see this, note that $26\cdot 15 + 51 = 441 = a_0^2$ and 
			$26\cdot 35 + 51 = 961 = a_1^2$ and suppose that for $k=n-1, n$, there exists 
			an $\eta \in \N$ such that $26\eta +51=a_k^2$.  Then, consider 
			$26\eta+51 = a_{n+1}^2$.  By definition $a_{n+1} = a_{n-1}+26$, therefore
			\[
				26\eta+51 = a_{n-1}^2+52a_{n-1} +26^2
			\]
			Solving for $a_{n-1}^2$ gives
			\[
				26(\eta-2a_{n-1}-26)+51 = a_{n-1}^2.
			\] 
			By induction, this has a solution in $\N$, call it $\ell$.  Then, 
			we can take $\eta = \ell +2a_{n-1}+26$.  It follows that each $a_n$ in 
			$\left\{a_n\right\}_{n=0}^\infty$ admits a distinct positive integer valued solution to 
			$26\eta+51=a_n^2$ and thus $M(0;[m_1+1,n_1+2],[m_2+1,n_2+2],q)$ 
			(with the restrictions outlined in Theorem \ref{main}) contains an infinite 
			subfamily of knots whose signatures are 0 and whose determinants are square.
		\end{proof}

\bibliographystyle{alphanum}
\bibliography{biblio}

\end{document}